\documentclass[a4paper,12pt,reqno]{amsart}
\usepackage{amsmath,amsthm,amssymb}
\usepackage{graphicx, color}
\usepackage{url}
\usepackage[numbers,sort&compress]{natbib}
\nonstopmode \numberwithin{equation}{section}

\allowdisplaybreaks

\allowdisplaybreaks

\begin{document}

\title[{Some summation theorems for Clausen's hypergeometric functions }]{Some summation theorems for Clausen's hypergeometric functions with unit argument}

\author[M.I. Qureshi and M. Shadab]{Mohammad Idris Qureshi and Mohd Shadab$^{*}$}

\address{Mohammad Idris Qureshi:   Department of Applied Sciences and Humanities,
 Faculty of Engineering and Technology,
 Jamia Millia Islamia (A Central University),
 New Delhi 110025, India}
\email{miqureshi\_delhi@yahoo.co.in}

\address{Mohd Shadab: Department of Applied Sciences and Humanities,
 Faculty of Engineering and Technology,
 Jamia Millia Islamia (A Central University),
 New Delhi 110025, India}
\email{shadabmohd786@gmail.com}



\thanks{*Corresponding author}

\begin{abstract}

Motivated by the work on hypergeometric summation theorems (recorded in the table III of Prudnikov et al. pp. 541-546), we have established some new summation theorems for Clausen's hypergeometric functions with unit argument in terms of $\pi$ and natural logarithms of some rational and irrational numbers. Also, we have given some modified summation theorems for Clausen's hypergeometric functions.

 \end{abstract}

\maketitle

\section{Introduction}\label{Intro}

The aim of this research article is to establish an elegant and manifestly relevant summation formula for the
Clausen’s series (see, e. g. \cite{Bailey1,Rainville,Slater,Sri4})

\begin{eqnarray}\label{eq(1.1)}
{_3}F_{2}\left[\begin{array}{r} 1, 1, \frac{p+q}{q};\\~\\2, \frac{p+2q}{q} ;\end{array}\  1\right]
\end{eqnarray}

with unit argument  and to support interest in generalized hypergeometric functions of this type.\\

Indeed, the classical field of hypergeometric functions ${_pF_q(z)}$ has been recently achieved a substantial progress by investigating the generalizing and producing various relationships between them. Often, the studied functions contain the integers and fractions in their numerator and denominator parameters in different ways, see mentioned papers, \cite{Gradshteyn,Krattenhalter,Lavoi1,Milgram1,Milgram2,Milgram3,Milgram4,Miller3,Miller1,Miller2,Prudnikov,Rainville,Rao,Shpot}( see, also references there in) .\\
\vskip.2cm

With reference to the importance of applications in mathematics, statistics, and mathematical physics, the hypergeometric functions readily reduce to a gamma function. On the basis of applications of hypergeometric functions have already been explored by a large number of authors notably C.F. Gauss, E.E. Kummer, S. Ramanujan, and others.

\vskip.2cm

 The papers \cite{Gradshteyn,Karl,Kim1,Krattenhalter,Lavoi1,Lewanowicz1,Rathie3,Rathie2,Rathie1} mentioned at the beginning, as well as the present one, discuss the summation formulas for the
function ${_3F_2(1)}$ that belong to the category of Clausen's hypergeometric function.
\vskip.2cm

We are now motivated enough by the work in the directions of summation theorems for Clausen's hypergeometric functions indicated in \cite{Milgram1,Milgram2,Milgram3,Milgram4,Miller3}. Term by term integration of special polynomials and Appell functions appear in such expansions in \cite{Appell,Chu1,Gauss}. Owing to numerous potential applications both in theoretical physics and mathematics, we believe that function \eqref{eq(1.1)}
or related deserves to be studied in the best way.

\section{Preliminaries}
The widely-used Pochhammer symbol $(\lambda)_{\nu}$ ~$(\lambda, \nu \in\mathbb{C})$ is defined by
\begin{equation}
\left(\lambda\right)_{\nu}:=\frac{\Gamma\left(\lambda+\nu\right)}{\Gamma\left(\lambda\right)}=\begin{cases}
\begin{array}{c}
1\\
~\\
\lambda\left(\lambda+1\right)\ldots\left(\lambda+n-1\right)
\end{array} & \begin{array}{c}
\left(\nu=0;\lambda\in\mathbb{C}\setminus\left\{ 0\right\} \right)\\
~\\
\left(\nu=n\in\mathbb{N};\lambda\in\mathbb{C}\right),
\end{array}\end{cases}
\end{equation}
it being understood $conventionally$ that $\left(0\right)_{0}=1$, and assumed $tacitly$ that the $\Gamma$ quotient exists.

The {\it generalized hypergeometric function} ${_p}F_q$, is defined by
\begin{eqnarray}\label{eq.(1.11)}
{_p}F_{q} \left[ \begin{array}{r}(a_p); \\ (b_q); \end{array} z \right] = \sum_{m=0}^{\infty} \frac{[(a_p)]_m}{[(b_q)]_m} \frac{z^m}{m!},
\end{eqnarray}

\begin{itemize}
  \item p and q are positive integers or zero,
  \item $z$ is a complex variable,
  \item $(a_p)$ designates the set ${a_1,a_2, . . . , a_p}$,
  \item the numerator parameters $a_1, . . . , a_p \in \mathbb{C}$ and the denominator parameters $b_1, . . . , b_q \in \mathbb{C}\setminus {\mathbb{Z}^{-}_{0}}, $
  \item $[(a_r)]_k = \displaystyle\prod_{i=1}^{r} (a_i)_k$. By convention, a product over the empty set is 1,
  \item $(a)_k$ is the Pochhammer's symbol.
\end{itemize}

Thus, if a numerator parameter is a negative integer or zero, the $_pF_q$ series terminates, and then we are led to a generalized hypergeometric polynomial.\\

\vskip.2cm
In 1856, Karl Weierstrass gave a different definition of gamma function
\begin{eqnarray}\label{eq(1.2)}
\frac{1}{\Gamma{(z)}}&=&z \exp{(\gamma z)} \prod_{n=1}^{\infty}\left[\left(1+\frac{z}{n}\right)\exp{\left(-\frac{z}{n}\right)}\right],
\end{eqnarray}
where $\gamma= 0.577215664901532860606512090082402431042\dots$ is called Euler-Mascheroni constant, and
$\frac{1}{\Gamma{(z)}}$ is an entire function of $z$, and
\begin{eqnarray}
\gamma&=&\lim_{n\rightarrow\infty}\left(1+\frac{1}{2}+\frac{1}{3}+.....+\frac{1}{n}-\ell n\,{(n)}\right)\nonumber.
\end{eqnarray}

The function

\begin{eqnarray}
\psi(z)=\frac{d}{dz}\{\ell n \,{\Gamma(z)}\}=\frac{\Gamma ^{\prime}(z)}{\Gamma(z)},
\end{eqnarray}
or, equivalently
\begin{eqnarray}
\ell n \,{\Gamma(z)}=\int_{1}^{z}\psi(\zeta)d\zeta,
\end{eqnarray}

is the logarithmic derivative of the gamma function (Psi function or digamma function).\\

\vskip.2cm
In 1813, Gauss \cite{Gauss} discovered an interesting formula for digamma (Psi) function as follows

\begin{eqnarray}\label{eq(Gauss)}
\psi(p/q) = -\gamma - \ell n\,{(q)}-\frac{\pi}{2}\cot{\left(\frac{\pi p}{q}\right)}+\sum_{j=1}^{[\frac{q}{2}]}{'}\left\{\cos\left({\frac{2\pi j p}{q}}\right)\ell n\, {\left(2-2\cos{\frac{2\pi j}{q}}\right)}\right\},
 \end{eqnarray}

where $1 \le p < q$ and $p, q$ are positive integers, and accent(prime) to right of the summation sign indicates the term corresponding to (last term) $j = \frac{q}{2}$ (when $q$ is positive even integer) should be divided by 2.\\
\vskip.2cm

In 2007, a simplified treatment of the Gauss formula was made by  Murty and Saradha \cite[p. 300, after eq.(4)]{Murty} (see also, Lehmer \cite[p. 135, after eq.(20)]{Lehmer} ) such that

\begin{eqnarray}\label{eq(Murty)}
\psi(p/q)=-\gamma-\ell n\,{(2q)}-\frac{\pi}{2}\cot{\left(\frac{\pi p}{q}\right)}+2\sum_{j=1}^{[\frac{q}{2}]}\left\{\cos{\left(\frac{2\pi pj}{q}\right)}\ell n\,{ \sin{\left(\frac{\pi j}{q} \right)}}\right\},
\end{eqnarray}
where $p = 1, 2, 3, \dots ,(q-1) ,\,  q = 2, 3, 4, \dots ; (p,q)= 1$.\\

\vskip.2cm

In volume III of Prudnikov et al. \cite[pp. 541-546]{Prudnikov}, summation theorems for Clausen's hypergeometric functions

\begin{eqnarray*}
&&{_3F_2} \left[1,1,\frac{1}{4}; 2,\frac{5}{4};1\right],\,\, {_3F_2} \left[1,1,\frac{1}{3}; 2,\frac{4}{3};1\right],\,\,{_3F_2} \left[1,1,\frac{3}{8}; 2,\frac{11}{8};1\right], \,\,{_3F_2} \left[1,1,\frac{1}{2}; 2,\frac{3}{2};1\right],\\
&&{_3F_2} \left[1,1,\frac{5}{8}; 2,\frac{13}{8};1\right], \,\,{_3F_2} \left[1,1,\frac{2}{3}; 2,\frac{5}{3};1\right], \,\,{_3F_2} \left[1,1,\frac{3}{4}; 2,\frac{7}{4};1\right], \,\,{_3F_2}\left[1,1,\frac{9}{8}; 2,\frac{17}{8};1\right], \\
&&{_3F_2}\left[1,1,\frac{7}{8}; 2,\frac{15}{8};1\right], \,\,{_3F_2} \left[1,1,\frac{5}{2}; 2,\frac{7}{2};1\right], \,\,{_3F_2}\left [1,1,\frac{5}{4}; 2,\frac{9}{4};1\right], \,\,{_3F_2} \left[1,1,\frac{3}{2}; 2,\frac{5}{2};1\right], \\
&&{_3F_2} \left[1,1,\frac{11}{8}; 2,\frac{19}{8};1\right], \,\,{_3F_2} \left[1,1,\frac{7}{4}; 2,\frac{11}{4};1\right],
\end{eqnarray*}
 are available.

\vskip.2cm

In the paper \cite{Qureshi}, we have given  summation theorems for Clausen's hypergeometric functions
${_3F_2} [1,1,n; 2,n+1;1]$, where n = 3,4,5,\dots,52.

\vskip.2cm

Motivated by the work, recorded in the table of Prudnikov et al. \cite[pp. 541-546]{Prudnikov} and a paper of Qureshi et al. \cite{Qureshi}, we have given some new and modified summation theorems for Clausen's hypergeometric functions using a celebrated formula of Gauss for digamma function in sections 4 and 5 respectively.

\section{Main Result}

In this section, first of all, we establish an interesting formula \eqref{eq(3.4)} for digamma function in the form of hypergeometric function, and connecting it with Gauss formula \eqref{eq(Murty)} for digamma function, we obtain an interesting and new result \eqref{eq(main-2.5)}, which plays a key role in our investigation.
\vskip.2cm

Let us recall the Weierstrass definition of gamma function \eqref{eq(1.2)}

\begin{eqnarray}
\frac{1}{\Gamma{(z)}}&=&z \exp{(\gamma z)} \prod_{n=1}^{\infty}\left[\left(1+\frac{z}{n}\right)\exp{\left(-\frac{z}{n}\right)}\right]\nonumber
\end{eqnarray}

\text{ Takin logarithm to the base e, we get }
\begin{eqnarray}\label{eq.(2.2)}
-\ell n \,{\Gamma{(z)}}&=&\ell n \,{z}+\gamma z+\sum_{n=1}^{\infty}\left\{\ell n \,\left(1+\frac{z}{n}\right)-\frac{z}{n}\right\}
\end{eqnarray}
Now, differentiating with respect to 'z', we get
\begin{eqnarray}
\frac{\Gamma^{'}{(z)}}{\Gamma{(z)}}&=&-\frac{1}{z}-\gamma -\sum_{n=1}^{\infty}\left\{\frac{1}{n+z}-\frac{1}{n}\right\}\nonumber\\
\frac{\Gamma^{'}{(z)}}{\Gamma{(z)}}&=&-\frac{1}{z}-\gamma +\sum_{n=1}^{\infty}\left\{\frac{z}{n(n+z)}\right\}
\end{eqnarray}
Using the definition \eqref{eq.(1.11)}, we get
\begin{eqnarray}
\psi(z)=\frac{\Gamma^{'}{(z)}}{\Gamma{(z)}}&=&-\frac{1}{z}-\gamma +\left(\frac{z}{1+z}\right){_3}F_{2}\left[\begin{array}{r} 1, 1, 1+z;\\
~\\
2, 2+z ;\end{array}\  1\right]
\end{eqnarray}
 where $z\ne 0, -1, -2, -3,\dots$ and $\psi(z)$ denotes the Psi (or Digamma) function.\\
\vskip.4cm

In eq(3.3) put $z= \frac{p}{q}$, we get
\begin{eqnarray}\label{eq(3.4)}
\psi\left(\frac{p}{q}\right) &=& - \gamma - \frac{q}{p} + \left(\frac{p}{p+q}\right)\,\,{_3}F_{2}\left[\begin{array}{r} 1, 1, \frac{p+q}{q};\\
~\\
2, \frac{p+2q}{q} ;\end{array}\  1\right]
\end{eqnarray}

\vskip.2cm
 Comparing \eqref{eq(Murty)} and \eqref{eq(3.4)}, we get

\begin{eqnarray}\label{eq(main-2.5)}
{_3}F_{2}\left[\begin{array}{r} 1, 1, \frac{p+q}{q};\\~\\2, \frac{p+2q}{q} ;\end{array}\  1\right]
= \left( \frac{q+p}{p}\right)\left[\frac{q}{p}-\ell n \,{2q}-\frac{\pi}{2}\cot{\left(\frac{\pi p}{q}\right)}+2\sum_{j=1}^{[\frac{q}{2}]}\left\{\cos{\left(\frac{2\pi pj}{q}\right)}\ell n \,{ \sin{\left(\frac{\pi j}{q} \right)}}\right\} \right]\nonumber\\,
\end{eqnarray}
$1 \le p < q.$

\section{Some new summation theorems for Clausen's function }

We have derived some new summation theorems for Clausen's hypergeometric function.

\begin{eqnarray}
&{_3}F_{2}&\left[\begin{array}{r} 1, 1, \frac{-1}{2};\\
~\\
2, \frac{1}{2} ;\end{array}\  1\right] =  \left\{  \frac{2}{3} -  \frac{2}{3}\ell n \,2 \right\}\\
&{_3}F_{2}&\left[\begin{array}{r} 1, 1, \frac{7}{2};\\
~\\
2, \frac{9}{2} ;\end{array}\  1\right] =  \frac{7}{5} \left\{  \frac{46}{15} -  2\ell n \,2 \right\}\\
&{_3}F_{2}&\left[\begin{array}{r} 1, 1, -\frac{4}{3};\\
~\\
2, -\frac{1}{3} ;\end{array}\  1\right] = \frac{4}{7} \left\{ \frac{15}{4} + \frac{\pi\sqrt{3}}{6} - \frac{3}{2}\ell n \,3 \right\}\\
&{_3}F_{2}&\left[\begin{array}{r} 1, 1, \frac{4}{3};\\
~\\
2, \frac{7}{3} ;\end{array}\  1\right] =  \left\{ 12 - \frac{2\pi}{\sqrt{3}} - 6\ell n \,3 \right\}\\
&{_3}F_{2}&\left[\begin{array}{r} 1, 1, \frac{10}{3};\\
~\\
2, \frac{13}{3} ;\end{array}\  1\right] = \frac{10}{7} \left\{ \frac{117}{28} - \frac{\sqrt{3}\pi}{6} - \frac{3}{2}\ell n \,3 \right\}\\
&{_3}F_{2}&\left[\begin{array}{r} 1, 1, \frac{6}{5};\\
~\\
2, \frac{11}{5} ;\end{array}\  1\right] = 6 \left\{ 5- \ell n \,10 - \left(\frac{1+\sqrt{5}}{\sqrt{(10-2\sqrt{5})}}\right)\frac{\pi}{2} + \frac{1}{2}\left( \sqrt{5} \ell n \,\left(\frac{\sqrt{5} - 1}{2}\right) -    \ell n \,\frac{\sqrt{5}}{4}  \right)     \right\}\nonumber\\
~~~~~~~~~~~~~~~\\
&{_3}F_{2}&\left[\begin{array}{r} 1, 1, \frac{7}{5};\\
~\\
2, \frac{12}{5} ;\end{array}\  1\right] = \frac{7}{2} \left\{ \frac{5}{2}- \ell n \,10 - \left(\frac{\sqrt{5} - 1}{\sqrt{(10+2\sqrt{5})}} \right)\frac{\pi}{2} + \frac{1}{2} \left( \sqrt{5} \ell n \,\left(\frac{\sqrt{5} + 1}{2} \right)-    \ell n \,\frac{ \sqrt{5} } {4}  \right)     \right\}\nonumber\\
~~~~~~~~~~~~~~~\\
&{_3}F_{2}&\left[\begin{array}{r} 1, 1, \frac{8}{5};\\
~\\
2, \frac{13}{5} ;\end{array}\  1\right] = \frac{8}{3} \left\{ \frac{5}{3} - \ell n \,10 + \left(\frac{\sqrt{5} - 1}{\sqrt{(10+2\sqrt{5})}} \right)\frac{\pi}{2} + \frac{1}{2} \left( \sqrt{5} \ell n \,\left(\frac{\sqrt{5} + 1}{2}\right) -    \ell n \,\frac{ \sqrt{5} } {4}  \right)     \right\}\nonumber\\
~~~~~~~~~~~~~~~\\
&{_3}F_{2}&\left[\begin{array}{r} 1, 1, \frac{9}{5};\\
~\\
2, \frac{14}{5} ;\end{array}\  1\right] = \frac{9}{4} \left\{ \frac{5}{4} - \ell n \,10 + \left(\frac{\sqrt{5} + 1}{\sqrt{(10-2\sqrt{5})}} \right)\frac{\pi}{2} + \frac{1}{2} \left( \sqrt{5} \ell n \,\left(\frac{\sqrt{5} - 1}{2} \right)-    \ell n \,\frac{ \sqrt{5} } {4}  \right)     \right\}\nonumber\\
~~~~~~~~~~~~~~~\\
&{_3}F_{2}&\left[\begin{array}{r} 1, 1, \frac{1}{6};\\
~\\
2, \frac{7}{6} ;\end{array}\  1\right] =  \left\{ \sqrt{3}\frac{\pi}{10} + \frac{3}{10}\ell n \,3 + \frac{2}{5}\ell n \,2 \right\}\\
&{_3}F_{2}&\left[\begin{array}{r} 1, 1, \frac{7}{6};\\
~\\
2, \frac{13}{6} ;\end{array}\  1\right] = 7 \left\{ 6 - \ell n \,12 - \sqrt{3}\frac{\pi}{2} - \ell n \,\sqrt{3} \right\}\\
&{_3}F_{2}&\left[\begin{array}{r} 1, 1, \frac{11}{6};\\
~\\
2, \frac{17}{6} ;\end{array}\  1\right] = \frac{11}{5} \left\{ \frac{6}{5} - \ell n \,12 + \sqrt{3}\frac{\pi}{2} - \ell n \,\sqrt{3} \right\}\\
&{_3}F_{2}&\left[\begin{array}{r} 1, 1, \frac{11}{10};\\
~\\
2, \frac{21}{10} ;\end{array}\  1\right] = 11 \left\{ 10 - \ell n \,20 -  \left(\frac{\sqrt{(10+2\sqrt{5})}}{\sqrt{5} - 1} \right)\frac{\pi}{2} + \frac{1}{2} \left( \sqrt{5} \ell n \,(\sqrt{5} - 2) -   \ell n \,\sqrt{5}  \right)     \right\}\nonumber\\
~~~~~~~~~~\\
&{_3}F_{2}&\left[\begin{array}{r} 1, 1, \frac{13}{10};\\
~\\
2, \frac{23}{10} ;\end{array}\  1\right] = \frac{13}{3} \left\{ \frac{10}{3} - \ell n \,20 -  \left(\frac{\sqrt{(10-2\sqrt{5})}}{\sqrt{5} + 1} \right)\frac{\pi}{2} + \frac{1}{2} \left( \sqrt{5} \ell n \,(\sqrt{5} + 2) -    \ell n \,\sqrt{5}  \right)     \right\}\nonumber\\
~~~~~~~~~~\\
&{_3}F_{2}&\left[\begin{array}{r} 1, 1, \frac{17}{10};\\
~\\
2, \frac{27}{10} ;\end{array}\  1\right] = \frac{17}{7} \left\{ \frac{10}{7} - \ell n \,20 +  \left(\frac{\sqrt{(10-2\sqrt{5})}}{\sqrt{5} + 1} \right)\frac{\pi}{2} + \frac{1}{2} \left( \sqrt{5} \ell n \,(\sqrt{5} + 2) -   \ell n \,\sqrt{5}  \right)     \right\}\nonumber\\
~~~~~~~~~~\\
&{_3}F_{2}&\left[\begin{array}{r} 1, 1, \frac{19}{10};\\
~\\
2, \frac{29}{10} ;\end{array}\  1\right] = \frac{19}{9} \left\{ \frac{10}{9} - \ell n \,20 +  \left(\frac{\sqrt{(10+2\sqrt{5})}}{\sqrt{5} - 1} \right)\frac{\pi}{2} + \frac{1}{2} \left( \sqrt{5} \ell n \,(\sqrt{5} - 2) -    \ell n \,\sqrt{5}  \right)     \right\}\nonumber\\
~~~~~~~~~~\\
&{_3}F_{2}&\left[\begin{array}{r} 1, 1, \frac{13}{12};\\
~\\
2, \frac{25}{12} ;\end{array}\  1\right] = 13 \left\{ 12 - \ell n \,24 -  \left(2+\sqrt{3}  \right)\frac{\pi}{2} +  \sqrt{3} \ell n \,(2-\sqrt{3}) -    \ell n \,\sqrt{3}    \right\}\\
&{_3}F_{2}&\left[\begin{array}{r} 1, 1, \frac{17}{12};\\
~\\
2, \frac{29}{12} ;\end{array}\  1\right] = \frac{17}{5} \left\{ \frac{12}{5} - \ell n \,24 -  \left(2 - \sqrt{3}  \right)\frac{\pi}{2} +  \sqrt{3} \ell n \,(2+\sqrt{3}) -    \ell n \,\sqrt{3}    \right\}\\
&{_3}F_{2}&\left[\begin{array}{r} 1, 1, \frac{19}{12};\\
~\\
2, \frac{31}{12} ;\end{array}\  1\right] = \frac{19}{7} \left\{ \frac{12}{7} - \ell n \,24 +  \left(2 - \sqrt{3}  \right)\frac{\pi}{2} +  \sqrt{3} \ell n \,(2+\sqrt{3}) -   \ell n \,\sqrt{3}    \right\}\\
&{_3}F_{2}&\left[\begin{array}{r} 1, 1, \frac{23}{12};\\
~\\
2, \frac{35}{12} ;\end{array}\  1\right] = \frac{23}{11} \left\{ \frac{12}{11} -\ell n \,24 +  \left(2 + \sqrt{3}  \right)\frac{\pi}{2} +  \sqrt{3} \ell n \,(2-\sqrt{3}) -    \ell n \,\sqrt{3}    \right\}
\end{eqnarray}

{\bf Independent proof of summation theorem (4.13):}\\
On substituting $z=\frac{1}{10}$ in eq. (3.3), we get

\begin{eqnarray}\label{eq.(4.1)}
\psi\left(\frac{1}{10} \right)= -\gamma - 10 +\frac{1}{11}{_3}F_{2}\,\left[\begin{array}{r} 1, 1, \frac{11}{10};\\
~\\
2, \frac{21}{10} ;\end{array}\  1\right].
\end{eqnarray}

Now, we calculate the value of $\psi\left(\frac{1}{10} \right)$ using the formula \eqref{eq(Murty)} as follows
\begin{eqnarray}\label{eq.(4.2)}
\psi\left(\frac{1}{10} \right ) = -\gamma -\ell n \, 20 - \left(\frac{\sqrt{(10+2\sqrt{5})}}{\sqrt{5} -1} \right) \frac{\pi}{2} +  \frac{1}{2}\{ \sqrt{5} \ell n \, (\sqrt{5} - 2) - \ell n \, \sqrt{5} \}.
\end{eqnarray}

On comparing the equations \eqref{eq.(4.1)} and \eqref{eq.(4.2)}, we get the desired result (4.13).\\

Similarly, we can derive other summation theorems.

\section{Some modified summation theorems for Clausen's function }

In this section, we have presented some corrected forms of erroneous summation theorems given in the table of Prudnikov et al. \cite[pp. 541-546]{Prudnikov}.

\begin{eqnarray}
&&{_3}F_{2}\left[\begin{array}{r} 1, 1, \frac{13}{8};\\
~\\
2, \frac{21}{8} ;\end{array}\  1\right] \circeq \frac{13}{50} \left\{ 16 - 5(\sqrt{2} - 1)\pi - 40 \ell n \, 2 + 10\sqrt{2} \ell n \,(1+\sqrt{2}) \right\},\\
&&{_3}F_{2}\left[\begin{array}{r} 1, 1, \frac{15}{8};\\
~\\
2, \frac{23}{8} ;\end{array}\  1\right] \circeq \frac{15}{196} \left\{ 32 + 7(1+\sqrt{2})\pi - 112 \ell n \, 2 - 28\sqrt{2} \ell n \,(1+\sqrt{2}) \right\},\\
&&{_3}F_{2}\left[\begin{array}{r} 1, 1, \frac{5}{3};\\
~\\
2, \frac{8}{3} ;\end{array}\  1\right] \circeq \frac{5\sqrt{3}}{12} \left\{ 3\sqrt{3} (1-\ell n \,3) - \pi \right\},
\end{eqnarray}
where the symbol $\circeq$ exhibits the fact that each of the equations (5.1), (5.2) and (5.3) does not hold true as stated.\\

The modified forms of the above summation theorems are as follows

\begin{eqnarray}
&&{_3}F_{2}\left[\begin{array}{r} 1, 1, \frac{13}{8};\\
~\\
2, \frac{21}{8} ;\end{array}\  1\right] = \frac{13}{50} \left\{ 16 + 5(\sqrt{2} - 1)\pi - 40 \ell n \, 2 + 10\sqrt{2} \ell n \,(1+\sqrt{2}) \right\},\\
&&{_3}F_{2}\left[\begin{array}{r} 1, 1, \frac{15}{8};\\
~\\
2, \frac{23}{8} ;\end{array}\  1\right] = \frac{15}{196} \left\{ 32 + 14(1+\sqrt{2})\pi - 112 \ell n \, 2 - 28\sqrt{2} \ell n \,(1+\sqrt{2}) \right\}\\
&&{_3}F_{2}\left[\begin{array}{r} 1, 1, \frac{5}{3};\\
~\\
2, \frac{8}{3} ;\end{array}\  1\right] = \frac{5\sqrt{3}}{12} \left\{ 3\sqrt{3} (1-\ell n \,3) + \pi \right\},.
\end{eqnarray}

\vskip.2cm
{\bf Remark:} We can not obtain the summation theorems of sections 4 and 5 with the help of summation theorems of Whipple, Dixon,  Watson and other associated contiguous relations of ${_3F_2}$.

\vskip.2cm
{\bf Concluding Remark :}\\
We conclude our present investigation by observing that several hypergeometric summation theorems for Clausen's hypergeometric function with unit argument have been deduced using celebrated formulas \eqref{eq(Murty)} and \eqref{eq(3.4)} in an analogous manner.


\begin{thebibliography}{99}

\bibitem{Appell}  Appell, P. and Kamp\'{e} de F\'{e}riet, J.; \textit{Fonctions hyperg\'{e}om\'{e}triques et hypersph\'{e}riques. Polynomes d'Hermite}, Gauthier-Villars, Paris, 1926.

\bibitem{Bailey1}   Bailey, W.N.; \textit{Generalised Hypergeometric Series}, Cambridge, England: Cambridge University Press, 1935.

\bibitem{Chu1} Chu, W.; Analytical formulae for extended ${_3F_2}$-series of Watson-Whipple-Dixon with two extra integer parameters,
\textit{Math. Comp.}, {\bf 81(277)} (2012), 467-479.

\bibitem{Gauss} Gauss, C.F.; Disquisitiones generales circa seriem infinitam etc., \textit{Comm. Soc. reg. Sci. Gott. rec., Vol II}, (1813) pp. 1-46.; reprinted in Werke {\bf 3} (1866).

\bibitem{Gradshteyn} Gradshteyn, I.S. and Ryzhik, I.M.; \textit{ Table of integrals, series and products}, 8th ed., Academic Press Inc., San Diego, CA. 2014.

\bibitem{Karl} Karlsson, P.W.; Reduction of hypergeometric functions with integral parameter differences, \textit{Indagationes
Mathematicae (Proceedings)},  {\bf 77(3)} (1974), 195-198.

\bibitem{Kim1} Kim, Y.S. Rathie, A.K. and Paris, R.B.; On two Thomae-type transformations for hypergeometric series with
integral parameter differences, \textit{Math. Comm.}, {\bf 19} (2014), 111-118.

\bibitem{Krattenhalter} Krattenhalter, C. and Rivoal, T.; How can we escape Thomae’s relations, \textit{J. Math. Soc. Japan}, {\bf58(1)} (2006), 183-210.

\bibitem{Lavoi1}  Lavoi, J.L. Grondin, F. and Rathie, A.K.; Generalizations of Whipple’s theorem on the sum of a ${_3F_2}$, \textit{J. Comp. Appl. Math.}, {\bf 72(2)} (1996), 293-300.

\bibitem{Lehmer} Lehmer, D.H.; Euler constants for arithmetical progressions,\textit{ Acta Arith.},  {\bf 27} (1975), 125-142.

\bibitem{Lewanowicz1} Lewanowicz, S.; Generalized Watson’s summation formula for ${_3F_2(1)}$,\textit{ J. Comp. Appl. Math.}, {\bf 86(2)} (1997), 375-386.

\bibitem{Milgram1} Milgram, M.; On some sums of digamma and polygamma functions, \url{ arXiv:math/0406338 [math.CA]}.

\bibitem{Milgram2}  Milgram, M.; On hypergeometric ${_3F_2(1)}$, \url{a review, arXiv:1011.4546 [math.CA]}.

\bibitem{Milgram3} Milgram, M.; 447 instances of hypergeometric ${_3F_2(1)}$, \url{arXiv:1105.3126 [math.CA]}.

\bibitem{Milgram4} Milgram, M.; Comment on a paper of Rao et al., an entry of Ramanujan and a new ${_3F_2(1)}$, \textit{J. Comp. Appl. Math.},  {\bf 201(1)} (2007), 1-2.


\bibitem{Miller3} Miller, A.R. and Paris, R.B.; Certain transformations and summations for generalized hypergeometric
series with integral parameter differences, \textit{Int. Trans. Spec. Func.},  {\bf 22(1)} (2011), 67-77.


\bibitem{Miller1} Miller, A.R. and Paris, R.B.; Clausen’s series ${_3F_2(1)}$ with integral parameter differences and transformations
of the hypergeometric function ${_2F_2(x)}$, \textit{Int. Trans. Spec. Func.},  {\bf 23(1)} (2012), 21-33.



\bibitem{Miller2} Miller, A.R. and Srivastava, H.M.; Karlsson-Minton summation theorems for the generalized hypergeometric
series of unit argument, \textit{Int. Trans. Spec. Func.},  {\bf 21} (2010), 603-612.

\bibitem{Murty} Murty, M.R. and Saradha, N.; Transcendental values of the digamma function,\textit{ J. Num. Theo.}, {\bf 125} (2007), 298-318.

\bibitem{Prudnikov}  Prudnikov, A.P. Brychkov, Y.A. and Marichev, O.I.; \textit{ Integrals and Series, Vol. 3: More special functions}, Gordon and Breach Science Publishers, 1990.


\bibitem{Rainville}Rainville, E.D.; \textit{Special Functions,} The Macmillan Co. Inc.,New York,1960; Reprinted by Chelsea Publ. Co. Bronx, New York, 1971.

\bibitem{Rao} Rao, K.S. Berghe, G.V. and Krattenhalter, C.; An entry of Ramanujan on hypergeometric series in his notebooks,
\textit{J. Comp. Appl. Math.}, {\bf 173(2)} (2005), 239-246.

\bibitem{Rathie3} Rathie, A.K. and Paris, R.B.; Extension of some classical summation theorems for the generalized hypergeometric
series with integral parameter differences, \textit{J. Classical Anal.}, {\bf 3(2)} (2013), 109-127.


\bibitem{Rathie1} Rathie, A.K. and Paris, R.B.; A note on some summations due to Ramanujan, their generalization and some
allied series, \url{arXiv:1301.4359 [math.CV]}.

 \bibitem{Rathie2} Rathie, A.K. and Rakha, M.A.; A study of new hypergeometric transformations, \textit{J. Phys. A: Math. Theor.},
{\bf 41(44)} (2008), 445202.

\bibitem{Shpot} Shpot, M.A.; A massive Feynman integral and some reduction relations for Appell functions, \textit{J. Math. Phys.},
 {\bf 48(12)} (2007), 123512-1-13.

\bibitem{Slater} Slater, L.J.; \textit{Generalized Hypergeometric Functions}, Cambridge University Press, Cambridge, 1966.

\bibitem{Sri4}Srivastava, H.M. and  Manocha, H.L.; \textit{A Treatise on Generating Functions,} Halsted Press(Ellis Horwood Ltd., Chichester, U.K.) John Wiley and Sons, New York, Chichester, Brisbane and Toronto, 1984.

\bibitem{Qureshi} Qureshi, M.I. Jabee, S. and Shadab, M.; Truncated Gauss hypergeometric series and its application in digamma function, \textit{(Communicated)}.

\end{thebibliography}
\end{document}